\documentclass[twocolumn,10pt]{article}
\usepackage{latex8}
\usepackage{times}
\usepackage{epsf}
\usepackage{psfig}

\newcommand{\noi}{\noindent}
\begin{document}

\title{How to Evolve Safe Control Strategies}
\author{G. W. Greenwood\thanks{greenwood@ieee.org}~ and X. Song \\ Dept. of Electrical \&
Computer Engineering \\ Portland State University, Portland, OR
97207 \\ \\
appears in Proc. of the 2002 NASA/DOD Conference on
Evolvable Hardware, 129-130, 2002}
\date{}

\maketitle \thispagestyle{empty}

\bibliographystyle{latex8}

\begin{abstract}
Autonomous space vehicles need adaptive control strategies that
can accommodate unanticipated environmental conditions.  The
evaluation of new strategies can often be done only by actually
trying them out in the real physical environment. Consequently, a
candidate control strategy must be deemed safe---i.e., it won't
damage any systems---prior to being tested online.  
How to do this efficiently has been a challenging problem.

We propose using evolutionary programming in conjunction with a
formal verification technique (called model checking) to evolve
candidate control strategies that are guaranteed to be safe for
implementation and evaluation.

\end{abstract}

\Section{Introduction}

Control strategies are critical ingredients of a space mission
because they indicate what actions are to be taken by the
spacecraft in response to environmental conditions. Unfortunately,
control strategies defined at the beginning of a mission may have
to be modified later on.  The need for this modification may
be due to system failures that reduce functionality or because the
spacecraft has encountered unanticipated environmental conditions.

An appealing method for dealing with these undesirable situations
is to use a \emph{reconfigurable system}, which can adopt a
different functionality.  For example, a reconfigurable system
eliminates the need for redundant hardware---which consumes
precious space and weight---by simply modifying the existing
hardware to compensate for the failure. However, despite the
enormous advantages of reconfiguration, reconfiguration
information originating from Earth will probably not arrive in
time to do any good.

The real solution lies with \emph{adaptive systems}---i.e.,
systems capable of self-reconfiguration in response to faults or a
changing operational environment~\cite{ds1}. This adaption is
performed \textit{in-situ} (in place), thereby removing any
reliance on Earth-bound resources for reconfiguration information.

In this paper we propose a method for evolving new control strategies
in a way {guaranteed} {to}  {be} {safe} during
the reconfiguration process. Our method fully supports
\textit{in-situ} adaption of the strategies.

\Section{Discussion}



Control strategies can be evolved \emph{extrinsically}, where each
strategy is simulated, but only the best one is actually
implemented, or \emph{intrinsically}, where each candidate
strategy is downloaded into the system and exercised in the real
physical environment. \textit{In-situ} extrinsic evolution may be 
problematic because some closed-form objective function is
necessary to assess efficacy, but it may not always be possible to 
define a suitable one. Thus, in most cases intrinsic evolution may be the
only thing that makes sense. It is therefore absolutely essential
that the control strategy be 
safe---i.e., it does no harm to the
controller itself nor to any other system.  This safety check
must be made prior to 
testing the new strategy online.

Our approach is to evolve a series of deterministic \emph{finite
state machines} (FSMs), each encoding a potential new control
strategy.  \textit{Evolutionary programming} (EP)~\cite{porto}
is used to evolve these FSMs. The suitability of each strategy
will be assessed by actually trying it in the real physical
environment. However, only control strategies that pass a safety
check will be downloaded for evaluation. We will borrow {automatic
formal verification methods} to assess this safety. These methods
use mathematically provable techniques to characterize a system
without conducting exhaustive simulation or testing. Specifically,
we will rely on \emph{model checking} (MC)
techniques~\cite{burch92} to verify the safety of candidate FSMs
generated by EP.  Although model checking has been extensively
used in hardware design and software verification, to the best of
our knowledge no prior research effort in formal methods has
attempted the problem we consider here.

MC is a formal method that verifies if a system, modelled as a
FSM, adheres to a specified property.  The properties of interest
are encoded as temporal logic expressions, which expresses
properties that change over time \cite{mcmillan93}. There are many
different kinds of temporal logic but \emph{computation tree
logic} (CTL) is the most widely used with model checkers. The
basic idea is a safety property is expressed in ordinary Boolean
logic, and then special temporal operators are added for
describing future events.

\begin{figure}[htb]
\centerline{\psfig{figure=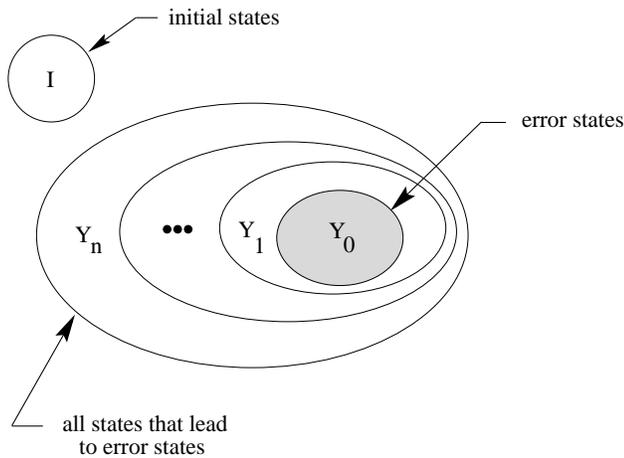,width=8.25cm,height=6.0cm}}
\caption{A graphical depiction of the model checking algorithm.
The control strategy is described by a FSM. $I$ is the set of all
FSM initial states and $Y_0$ is the set of FSM states that
violates a safety property. The algorithm recursively computes
$Y_{i+1} = \mbox{Pre}(Y_i)\bigcup~Y_i$ for $i=0,1,2,\ldots n-1$
where Pre($Y_i$) is the preimage of the set $Y_i$.  $Y_n$ then
represents the set of all FSM states that can reach an error
state. The system is safe if $Y_n\bigcap~I= \emptyset$. This check
can be done in linear time.} \label{mc}
\end{figure}

MC has been used to verify properties in systems with hundreds of
thousands of states. In practice, control strategies tend to have
orders of magnitude fewer states.  The MC algorithm complexity is
linear in the size of the FSM and in the length of the CTL
expression~\cite{burch94}, so the safety of a control strategy can
be quickly verified. A graphical representation of the MC
algorithm is shown in Figure \ref{mc}.

Several important issues concerning using MC to check safeness of
control strategies are worth highlighting:
\begin{itemize}
\item \textit{The large number of states in physical systems often
forces one to use a reduced FSM model where some details are
abstracted out. Model checking cannot guarantee safety under these
circumstances.}

In our approach the EP algorithm renders FSMs which are complete
in the sense that {\underline{every} aspect of the control
strategy is explicitly described in the FSM structure}. No details
are abstracted out or reduced so the safety check results are
guaranteed.

\item \textit{Model checkers typically provide trace information
to help pinpoint where the safety property failed.}

We will not use this feature.  In fact, we treat the entire safety
issue as a decision problem---i.e., either the strategy is safe or
it is not. Unsafe control strategies are immediately discarded, so
there is no need to know why it is unsafe.

\item \textit{Model checkers are used to verify functional
specifications and other properties, e.g., liveness.}

In our approach model checking only verifies safety. All other
performance criteria are assessed by trying out the control
strategy in the physical environment.

\end{itemize}

\Section{Implementation Details} \label{adv}


 \noi Our method can be summarized as follows:

\begin{itemize}
\item Control strategies are encoded with FSMs. \item An EP
algorithm generates candidate control strategies by evolving FSMs.
EP is ideally suited for this task~\cite{porto}. \item Safety properties are encoded
as CTL expressions. \item A symbolic model checker accepts the FSM
and CTL expressions as input, and quickly checks to see if the
control strategy is safe. The correctness of the safety check is
guaranteed. \item Safe control strategies are evaluated in the
physical environment whereas unsafe strategies are discarded.
\item The EP algorithm runs a fixed number of generations or 
terminates sooner if a suitable control strategy is found. 
The best performing FSM is implemented as the new control strategy.
\end{itemize}

\bibliography{nasa}

\end{document}